\documentclass{amsproc} 
\usepackage{pst-plot}
\usepackage{pstricks}
\usepackage{mathrsfs}

\allowdisplaybreaks

\def\Bin{{\rm Bin}}

\def\E{\mathsf{E}}
\def\P{\mathsf{P}}
\def\var{\mathsf{var}}
\def\cov{\mathsf{cov}}

\def\b{\beta}
\def\a{\alpha} 

\newtheorem{theorem}{Theorem}[section]
\newtheorem{lemma}[theorem]{Lemma}
\newtheorem{corollary}[theorem]{Corollary}
\newtheorem{proposition}[theorem]{Proposition}

\theoremstyle{definition}

\theoremstyle{remark}
\newtheorem{remark}[theorem]{Remark}

\numberwithin{equation}{section}

\newcommand{\euler}[2]{\genfrac{ < }{ > }{0pt}{}{#1}{#2}}
\newcommand{\stir}[2]{\genfrac{ \{ }{ \} }{0pt}{}{#1}{#2}}

\newcommand\cF{\mathcal F}
\newcommand\cL{\mathcal L}
\newcommand\cT{\mathcal T}
\newcommand{\pgf}{probability generating function}
\newcommand{\Sx}{\ensuremath{\mathsf S}}
\newcommand{\W}{\ensuremath{\mathsf W}}
\newcommand{\kk}{_{k-1}}
\newcommand{\nn}{_{n-1}}
\newcommand\eqd{\overset{\mathrm{d}}{=}}
\newcommand\dto{\overset{\mathrm{d}}{\longto}}
\newcommand{\longto}{\longrightarrow}
\newcommand{\Hh}{h}

\newcommand\bigpar[1]{\bigl(#1\bigr)}
\newcommand\Bigpar[1]{\Bigl(#1\Bigr)}

\newcommand\lrpar[1]{\left(#1\right)}

\def\rompar(#1){\textup(#1\textup)}    

\newcommand{\refT}[1]{Theorem~\ref{#1}}
\newcommand{\refC}[1]{Corollary~\ref{#1}}
\newcommand{\refL}[1]{Lemma~\ref{#1}}

\newcommand{\datex}[1]{\thanks{\emph{Date}: #1}}

\begin{document}

\title{Asymptotic Normality of Statistics on  Permutation Tableaux}

\author{Pawe{\l} Hitczenko}
\address{Department of Mathematics, Drexel University, Philadelphia, Pennsylvania 19104}
\email{phitczenko@math.drexel.edu}
\thanks{The first author was supported in part by the NSA Grant \#H98230-09-1-0062.}

\author{Svante Janson}
\address{Department of Mathematics,
Uppsala University, Box 480, SE-751 06 Uppsala, Sweden}
\email{svante.janson@math.uu.se}
\thanks{Research of the second author was partially done while visiting Institut
  Mittag-Leffler, Djursholm, Sweden}

\datex{March 16, 2009}




\begin{abstract}
In this paper we use a probabilistic approach to derive the
expressions for the characteristic functions of basic statistics
defined on permutation tableaux. Since our expressions are exact, we
can identify the distributions of basic statistics (like the number of
unrestricted rows, the number of rows, and the number of 1s in the
first row) exactly. In all three cases the distributions are known to
be asymptotically normal after a suitable normalization.   We also
establish the asymptotic normality of the number of superfluous
1s. The latter results relies on a bijection between permutation
tableaux and permutations and on a rather general sufficient condition
for the central limit theorem for the sums of random variables in
terms of dependency graph of the summands.  
\end{abstract}

\maketitle

\section{Introduction}

Permutation tableaux are relatively new objects that are in bijection
with permutations \cite{Bu,CN,SW}. They were introduced in the context
of  enumeration of  the totally positive Grassmannian cells 
\cite{Postnikov, Williams}. More recently, permutation tableaux
generated additional research activity when they have been connected
in \cite{C,CW,CW1} 
to a particle model in statistical physics   called the Partially
ASymmetric Exclusion  Process (PASEP); see
\cite{BrakEssam,Corteel2,Derrida0,Derrida1,jumping,Sasamoto,Askey-Wilson}
for more information on PASEP. 
 
A {\em permutation tableau}  \cite{SW} is  a Ferrers diagram of a partition 
of a positive integer into non--negative parts
whose boxes are filled   with $0$s and
$1$s according to the following rules:
\begin{enumerate}
\item Each column of the diagram contains at least one $1$.
\item There is no $0$ which has a $1$ above it in the same column
{\em and} a $1$ to its left in the same row.
\end{enumerate}
An example is given in Figure \ref{fig}. 

\begin{figure}[h]
 \centering
        \psset{unit=0.4cm}
\begin{pspicture}(-1,0)(7,5)        
\psframe[dimen=middle](0,0)(1,1)
\psframe[dimen=middle](0,1)(1,2)
\psframe[dimen=middle](1,1)(2,2)
\psframe[dimen=middle](2,1)(3,2)
\psframe[dimen=middle](0,2)(1,3)
\psframe[dimen=middle](1,2)(2,3)
\psframe[dimen=middle](2,2)(3,3)
\psframe[dimen=middle](3,2)(4,3)
\psframe[dimen=middle](4,2)(5,3)
\psframe[dimen=middle](0,3)(1,4)
\psframe[dimen=middle](1,3)(2,4)
\psframe[dimen=middle](2,3)(3,4)
\psframe[dimen=middle](3,3)(4,4)
\psframe[dimen=middle](4,3)(5,4)
\psframe[dimen=middle](0,4)(1,5)
\psframe[dimen=middle](1,4)(2,5)
\psframe[dimen=middle](2,4)(3,5)
\psframe[dimen=middle](3,4)(4,5)
\psframe[dimen=middle](4,4)(5,5)
\psframe[dimen=middle](5,4)(6,5)
\psframe[dimen=middle](6,4)(7,5)
\rput(0.5,0.5){\small 1}
\rput(0.5,1.5){\small 0}
\rput(0.5,2.5){\small 0}
\rput(1.5,2.5){\small 1}
\rput(2.5,2.5){\small 1}
\rput(3.5,2.5){\small 1}
\rput(4.5,2.5){\small 1}
\rput(0.5,3.5){\small 0}
\rput(1.5,3.5){\small 0}
\rput(2.5,3.5){\small 1}
\rput(3.5,3.5){\small 0}
\rput(4.5,3.5){\small 1}
\rput(0.5,4.5){\small 0}
\rput(1.5,4.5){\small 0}
\rput(3.5,4.5){\small 0}
\rput(4.5,4.5){\small 0}
\rput(2.5,4.5){\small  1}
\rput(5.5,4.5){\small  1}
\rput(6.5,4.5){\small  1}
\rput(2.5,1.5){\small   0}
\rput(1.5,1.5){\small  0}
\psframe[dimen=middle](0,0)(0,-1)
\end{pspicture}
         \caption{Example of a permutation tableau}
 \label{fig}
\end{figure}

The size parameter of a permutation tableau is its {\em length}
defined as  the number of rows plus the number of columns of the
tableau.  
For example, the tableau in Figure \ref{fig} has 6 rows and 7 columns
so its length is 13.  

Different statistics on permutation tableaux were defined in \cite{CW1,SW}.
We recall that a 0 in a permutation tableau is  {\it restricted} if 
there is a 1 above it in the same column. A row is {\it unrestricted}
if it does not contain a restricted 0.  
A 1 is {\em superfluous} if it has a 1 above itself in the same column.
We will be interested in the number of unrestricted rows, the number
of superfluous 1s, as well as the number of 1s in the top row, and
the number of rows. For example the tableau in Figure \ref{fig} has   
three superfluous 1s, six rows, three 1s in the top row,  and five unrestricted rows.

Most of the past research on statistics of permutation tableaux was
based on bijections between permutation tableaux of length $n$ and
permutations of $[n]:=\{1,\dots,n\}$ and using known properties of
permutations. This was not always easy as sometimes it is not that
easy to see into what parameter a given statistic is mapped by a
bijection. In \cite{ch} a direct approach based on a probabilistic
consideration was proposed. It enabled the authors to compute the
expected values of  these statistics in a simple and unified way. In
this work we go one step further and we will  concentrate on the
limiting distributions of (properly normalized)  statistics mentioned
above. In most cases we will accomplish it by computing the
probability generating functions 
of the quantities in question. As a matter of fact, since
the approach proposed in \cite{ch} allows for exact (and not only
asymptotic) computation,  in most cases we will be able to identify
the distribution of a given statistics exactly and not only
asymptotically. As a consequence, we can see, in particular,  that:
\begin{itemize}
\item the number of unrestricted rows in a random permutation tableaux
  of length $n$ has the same distribution as the number of cycles, or
  the number of records, in a random permutation of $[n]$
  (a record in a permutation $\sigma=(\sigma_i)$ is any $\sigma_i$ such that
  $\sigma_i>\sigma_j$ for $j<i$);
\item the number of rows has the same distribution as the number of
  descents in a random permutation counted by Eulerian numbers (a
  descent in $(\sigma_i)$ is any pair $(\sigma_i,\sigma_{i+1})$ such
  that $\sigma_i>\sigma_{i+1}$);   
\item  the number of 1s in the first row is, in distribution,  one
  less than the number of unrestricted rows.  
\end{itemize}
The first two results follow from the analysis of known bijections
between permutations and permutation tableaux (see \cite{CN,SW}),
while the third is a consequence of an involution on permutation
tableaux presented in \cite{CW1}.   However, neither the involution
from \cite{CW1} nor any of the bijections from \cite{CN,SW} are very
straightforward. In addition, the first two results given above
required two different bijections in \cite{SW}. So, while these
bijections carry more information, they turn out to be rather
cumbersome to work with. Looked at from this perspective, our approach
provides a streamlined and unified derivation of the  above results.

As for the number of superfluous 1s, the situation is a bit more
complicated. Although we do derive its \pgf, 
we do not deduce  directly from this the asymptotic normality of its
distribution. Instead, in this case we will rely on the fact that one of
the two bijections between permutation tableaux and permutations
described in \cite{SW} sends the number of superfluous 1s to the
number of generalized patterns $31{-}2$ and we will prove the central
limit theorem for the number of such patterns in random permutation of
$[n]$ (an occurrence of a generalized pattern $31{-}2$ in $(\sigma_i)$
is any pair $(i,j), 1<i<j$ such that
$\sigma_{i-1}>\sigma_j>\sigma_i$). Our proof is based on a rather
general sufficient condition for the central limit theorem developed
in \cite{SJ58} and \cite[Section~6.1]{jlr}. 
B\'ona \cite{Bona} has used the same method to show similar results on
asymptotic normality of general permutation patterns, but
as far as we know the central limit
theorem for the number of a generalized pattern such as $31{-}2$ is
new. Furthermore, the number of permutation tableaux of length $n$
with $k$ superfluous 1s 
is equal to the number of permutations of $[n]$ with $k$ crossings 
(a crossing \cite{C} in a permutation 
$(\sigma_i)$ is a pair
$(i,j)$ such that $i<j\le \sigma_i<\sigma_j$ or
$i>j>\sigma_i>\sigma_j$). It follows, therefore, that the number of
crossings in a random permutation of $[n]$ is asymptotically normal.

\section{Basic facts}

Let $\mathcal T_n$ be the set of all permutation tableaux of length
$n$. 
We denote the uniform probability measure on $\mathcal T_n$ by $\P_n$
and $\E_n$ will denote the expectation with respect to $\P_n$. We
denote by $R_n$, $C_n$, $U_n$, $F_n$, and $S_n$ the random variables 
representing the numbers of rows, columns, 
unrestricted rows, 
1s in the top row, and superfluous 1s,
respectively, in a random tableau of length $n$. 
(We may allow $n=0$; there is a single empty permutation tableau of
length 0, and $U_0=F_0=R_0=C_0=S_0=0$. It is sometimes convenient to
start inductions with $n=0$, but for simplicity we often treat $n\ge1$ only.)

If $\cF$ is a  
$\sigma$--algebra  then $\E(\ \cdot\  |\cF)$ denotes the conditional
expectation 
given $\cF$.  If
$X$ is
non-negative integer valued, for example one of the permutation
tableaux statistics just defined,
we  let $g_X(z)=\E
z^X$ be its  probability generating function (in general defined at
least for $|z|\le1$; for the variables consider here,
$z$ can be any complex number).
We will often omit the subscript $X$.

The arguments in \cite{ch} were based on a construction of tableaux of
size $k$ from tableaux of size $k-1$ by extending the length of the
latter from its  
south-west (SW) corner either to the south (creating a new empty row) or to
the west (creating and then filling a new column). 
Note that each permutation tableau in $\cT_k$ is an extension of a
 unique permutation tableau in $\cT_{k-1}$. 
We refer the reader to \cite[Section~2]{ch} for details.
We let $M_k$ indicate the direction of the move at the  
 $k$th step (i.e. when the length is increased from $k-1$ to $k$). We
refer to $M_k$ as the $k$th move and will write $M_k=\Sx$ or $M_k=\W$ to
indicate its direction.  
 
The following  simple observations were crucial for the arguments
in \cite{ch} and are crucial here as well (see \cite{ch}
for some details that are omitted here). 
\begin{itemize}
\item
When extending a tableau by a \W{} move, the new column has to be filled
with 0s and 1s;
each restricted row must be filled with 0, but the unrestricted rows
can be filled arbitrarily except that there must be at least one 1.
Hence, a tableau with $U$ unrestricted rows can be extended by a
\W{} move in $2^U-1$ ways; since there always is a unique \Sx{} move, the
total number of extensions is $2^U$.
\item
We use
 $T_n$ to denote a generic element of $\cT_n$,
and use $T_{n-1}$ to denote the corresponding element
of $\cT_{n-1}$ (i.e., such that $T_n$ is an
 extension of $T_{n-1}$); we let
$U_n=U_n(T_n)$ be the number of unrestricted rows in $T_n$
 and 
$U_{n-1}=U_{n-1}(T_{n-1})$ the number of
 unrestricted rows in $T_{n-1}$, and similarly for the other variables
 that we study.
In this way,
 $\P_n$ induces a probability measure 
(also denoted $\P_n$) on $\mathcal T_{n-1}$, namely
 each element 
  of $\mathcal T_{n-1}$ is assigned a measure that is proportional to
  the number of tableaux from $\mathcal T_n$ it generates when its
  length is increased from $n-1$ to  $n$. 
Note that this differs from the uniform distribution
 $\P_{n-1}$ on $\cT_{n-1}$.
Since $|\mathcal T_n|=n!$ and a  tableau $T$ from $\mathcal
  T_{n-1}$  generates $2^{U_{n-1}(T)}$  tableaux of length $n$   
the relationship between these two measures is easy to find.
One way to state this relationship is that  if $X$ is any random
  variable on  $\mathcal T_{n-1}$ then  
\begin{equation}\label{relation}
\E_nX=\frac1{n}\E_{n-1}\left(
2^{U_{n-1}}X\right), \end{equation} 
where integration on the left is with respect to the measure induced
on $\mathcal T_{n-1}$ by $\P_n$ and integration on the right is with
respect to $\P_{n-1}$. 
Equivalently, if we let $\cF_{n-1}$ be the
$\sigma$-algebra on $\cT_n$ generated by the mapping $T_n\mapsto T_{n-1}$, then
$d\P_{n}/d\P_{n-1}=2^{U_{n-1}}/n$ on $\cF_{n-1}$.   
\item  
The sequence of the distributions of the number
of unrestricted rows is given by $U_0\equiv0$, 
 $ U_1\equiv1$ and for $n\ge1$ the conditional distribution of $U_{n}$
given $\cF_{n-1}$ (i.e., given $T_{n-1}$) is (under $\P_n$)
\begin{equation}\label{conddist}
\mathcal{L}(U_{n}\mid \cF_{n-1})=1+\Bin(U_{n-1}),\end{equation}
where $\Bin(m)$ denotes a binomial random variable with parameters $m$
and $1/2$. 
\end{itemize}
To illustrate how these  facts are put together to  work  consider
the number of unrestricted rows, $U_n$. For its \pgf{} we have:  
\begin{equation}\label{gun}
  \begin{split}
g_{U_n}(z)&=\E_nz^{U_n}=
\E_n\E_n\left( z^{U_n}\mid \cF_{n-1}\right)
=
 \E_n\E_n\left( z^{1+\Bin(U_{n-1})}\mid \cF_{n-1}\right) 
\\&=z\E_n\Bigpar{\frac{z+1}2}^{U_{n-1}}
= \frac{z}n\E_{n-1}\lrpar{2^{U_{n-1}}\Bigpar{\frac{z+1}2}^{U_{n-1}}}\\&=
\frac{z}n\E_{n-1}\left(z+1\right)^{U_{n-1}},	
  \end{split}
\end{equation}
where we have used (in that order) conditioning, (\ref{conddist}), the
obvious fact that 
$g_{\Bin(m)}(z)=\E z^{\Bin(m)}=\left(\frac{z+1}2\right)^m$, and
(\ref{relation}).  
It follows by induction, and $g_{U_1}(z)=z$, that
\begin{equation}
  \label{gunn}
g_{U_n}(z) =\frac{\Gamma(z+n)}{\Gamma(z)n!}
=\binom{z+n}n
=\prod_{j=0}^{n-1}\frac{z+j}{j+1}=\prod_{k=1}^n\Bigpar{1-\frac1k+\frac zk}.
\end{equation}
The factor on the right-hand side is the \pgf{} of a
random variable that is 1 with probability $1/k$ and 0 with
probability $1-1/k$. Since the product of \pgf{s}
corresponds to summing independent random variables, we obtain the
following statement 
\begin{theorem}\label{thm:U}
The number of unrestricted rows $U_n$ is distributed like
$$U_n\eqd\sum_{k=1}^nJ_k,$$ 
where $J_1$, $J_2,\dots$, are independent  indicators with $\P(J_k=1)=1/k$.
In particular, if $\Hh_n=\sum_{k=1}^nk^{-1}$ and 
$\Hh_n^{(2)}=\sum_{k=1}^nk^{-2}$ is the $n$th harmonic and generalized
harmonic number, respectively, then
$$\E_nU_n=\sum_{k=1}^n\P(J_k)=\Hh_n,
\quad\var_n(U_n)=\sum_{k=1}^n\var(J_k)
=\sum_{k=1}^n\frac1k\left(1-\frac1k\right)=\Hh_n-\Hh_n^{(2)},$$ 
and
$$\frac{U_n-\ln n}{\sqrt{\ln n}}\dto N(0,1),$$
where $N(0,1)$ denotes a standard normal random variable.
\end{theorem}
 
\begin{remark}
It is seen from the above statement that the  distribution of $U_n$
 coincides with the distribution of the number of cycles in a random
 permutation of $[n]$ (see for example \cite[Chapter X,
 Section~6(b)]{feller} or \cite[Chapter~4,
 Section~3]{cycles}),
or, equivalently,
the number of RL minima in a random permutation (which is known to
 be equidistributed with the number of cycles).   
Indeed, a bijection between 
 permutation tableaux and permutations  described in \cite{CN} maps
 unrestricted rows  onto RL minima in the corresponding permutation.
 \end{remark}

\begin{remark}
  More generally, we may define on $\cT_n$ permutation
  tableaux $T_k\in\cT_k$ for every $k\le n$, such
  that $T_{k+1}$ is an extension of $T_k$, and let
  $\cF_k$ be the $\sigma$-algebra generated by
  $T_k$. It can be seen by induction, similar to \eqref{gun}, 
that a permutation tableau $T_k$
  of length $k$ can be extended to $(n-k)!\,(n-k+1)^{U_k(T_k)}$
  permutation tableaux in $\cT_n$, and hence 
$d\P_n/d\P_k=(n-k+1)^{U_k}/\binom{n}{k}$  on $\cF_k$,   $1\le k\le  n$. 
Furthermore, \eqref{gunn} can be generalized to 
\begin{equation*}
\E_n(z^{U_n}\mid\cF_k)=\frac{\Gamma(z+n-k)k!}{\Gamma(z)n!}(z+n-k)^{U_k}.  
\end{equation*}
We can further study the sequence $(U_k)_{k=1}^n$ as a
stochastic process under $\P_n$; similar calculations show that this is an
inhomogeneous Markov process with transitions given by
$\cL(U_{k+1}\mid\cF_k)=1+\Bin(U_k,\frac{n-k}{n-k+1})$, which
generalizes \eqref{conddist}.
\end{remark}

\section{The number of 1s in the first row} 
Recall that $F_n$ denotes the number of 1s in the top row of a
random permutation tableau of length $n$. The following fact is a
consequence of an involution on permutation tableaux given in 
\cite{CW1}.  Since this involution is not too easy to describe we
provide a more direct justification of (\ref{U=F}) based on the
arguments given in the preceding section. 
\begin{theorem}\label{thm:fn} For every $n\ge 1$, we have
\begin{equation}\label{U=F}
F_n\eqd U_n-1.\end{equation}
In particular,
$$\E_nF_n=\Hh_n-1,
\qquad\var_n(F_n)=\Hh_n-\Hh_n^{(2)},$$ 
and
$$\frac{F_n-\ln n}{\sqrt{\ln n}}\dto N(0,1).$$
\end{theorem}

To establish (\ref{U=F}) it will be convenient to prove some auxiliary
lemmas. Define  $G_k$ to be  the position (counting only unrestricted
rows) of the topmost 1 on the $k$th move, provided that this move 
is \W. Then 
 $$F_n=\sum_{k=1}^nI_{G_k=1}.$$ Note that $G_k$ is undefined if the
$k$th move is \Sx. While it is inconsequential in this section as  we will
be interested in the event $\{G_k=1\}$, for the  
purpose of the subsequent sections it is convenient to set
$G_k=U_{k-1}+1$ if the $k$th move is \Sx. 
We have 
\begin{lemma}\label{lem:G} For all $k\ge2$, 
$$\P_k(G_k=1\mid \cF_{k-1})=\frac{1}{2}.$$
\end{lemma}
\begin{proof} 
A given tableau $T\in\cT_{k-1}$ has
$2^{U\kk(T)}$ extensions. The number of extensions by a
$\W$ move with a 1 in the topmost 
row (which always is unrestricted) is
$2^{U\kk(T)-1}$ since the other $U\kk(T)-1$
unrestricted rows can be filled arbitrarily. Hence the probability
that $G_k=1$ in a random extension of $T$ is $2^{U\kk(T)-1}/2^{U\kk(T)}=1/2$.
\end{proof}

\begin{lemma}\label{iter-F} 
For any complex numbers $z,w$ and every $k\ge2$, 
$$
\E_k(z^{I_{G_{k}}}w^{U_k}\mid \cF_{k-1})
=w\frac{z+w}{w+1}\left(\frac{w+1}2\right)^{U_{k-1}}.$$
\end{lemma}
\begin{proof}
Let $T\in \cT_{k-1}$. If $T$ is extended with
$G_k=1$, i.e., by a $\W$ move with 1 added in the
topmost  row, then we may
add 0 or 1 arbitrarily to all other unrestricted rows, and the rows
with 1 added remain unrestricted. Hence, conditionally given
$G_k=1$, the extension has $1+\Bin(U\kk(T)-1)$
unrestricted rows. 
Recall from \eqref{conddist} that without further
conditioning, the number of unrestricted rows is $1+\Bin(U\kk(T))$.
Consequently, using also Lemma \ref{lem:G},
\begin{align*}
\E_k(z^{I_{G_{k}}}w^{U_k}&\mid \cF_{k-1})=
z\E_k(w^{U_k}I_{G_{k}=1}\mid \cF_{k-1})
+\E_k(w^{U_k}I_{G_{k}\ne1}\mid \cF_{k-1})\\&=
(z-1)\E_k(w^{U_k}I_{G_{k}=1}\mid \cF_{k-1})
+\E_k(w^{U_k}\mid \cF_{k-1})
\\&=
(z-1)\E_k(w^{1+\Bin(U_{k-1}-1)}I_{G_{k}=1}\mid \cF_{k-1})
+\E_k(w^{1+\Bin(U_{k-1})}\mid \cF_{k-1})
\\&=
(z-1)w\left(\frac{w+1}2\right)^{U_{k-1}-1}\P_k(G_{k}=1\mid \cF_{k-1})
+w\left(\frac{w+1}2\right)^{U_{k-1}}
\\&=
(z-1)\frac w{w+1}\left(\frac{w+1}2\right)^{U_{k-1}}
+w\left(\frac{w+1}2\right)^{U_{k-1}}
\\&=
w\frac{z+w}{w+1}\left(\frac{w+1}2\right)^{U_{k-1}}.
\qedhere
\end{align*}
\end{proof}

\begin{lemma}
  \label{L-FU}
The joint \pgf{} of $F_n$ and $U_n$ is given by
\begin{equation*}
  \E_n\bigpar{z^{F_n}w^{U_n}}
=w\frac{\Gamma(z+w+n-1)}{\Gamma(z+w)n!}
=w\prod_{k=2}^n\frac{z+w+k-2}k.
\end{equation*}
\end{lemma}

\begin{proof}
By Lemma \ref{iter-F} and \eqref{relation},
  \begin{equation*}
	\begin{split}
\E_n\bigpar{z^{F_n}w^{U_n}}	
&=\E_n\E_n\bigpar{z^{F\nn+I_{G_n=1}}w^{U_n}\mid\cF\nn}
\\&
=\E_n\lrpar{z^{F\nn}w\frac{z+w}{w+1}\Bigpar{\frac{w+1}2}^{U\nn}}
\\&
=w\frac{z+w}{w+1}\frac1n\E\nn\lrpar{z^{F\nn}(w+1)^{U\nn}},
	\end{split}
  \end{equation*}
and the formula follows by induction.
\end{proof}

We now can complete the proof of Theorem~\ref{thm:fn}.

\begin{proof}[Proof of Theorem \ref{thm:fn}] 
Taking $w=1$ in Lemma \ref{L-FU} we obtain
\begin{equation*}
  \E_n z^{F_n}=\frac{\Gamma(z+n)}{\Gamma(z+1)n!},
\end{equation*}
which equals $\E_n z^{U_n-1}$ by \eqref{gunn}.
\end{proof}

Note that we recover \eqref{gunn} by taking $z=1$ in Lemma \ref{L-FU}.
We can also describe the joint distribution of $F_n$ and $U_n$.

\begin{theorem}
  \label{T-FU}
For every $n\ge1$, the joint distribution of $F_n$ and
$U_n$ is given by
\begin{equation*}
(F_n,U_n)\eqd\sum_{k=1}^n(J_k,I_k),  
\end{equation*}
where the random vectors $(J_k,I_k)$ are independent and $J_1=0$,
$I_1=1$, and, for $k\ge2$, 
$\P(I_k=1,\,J_k=0)=\P(I_k=0,\,J_k=1)=1/k$,
$\P(I_k=0,\,J_k=0)=1-2/k$.
\end{theorem}

\begin{proof}
  We have $\E(z^{J_k}w^{I_k})=(k-2+z+w)/k$ if $k\ge2$ and
  $w$ if $k=1$, and thus the joint \pgf{} of the
  right-hand side equals the product in Lemma \ref{L-FU}.
\end{proof}

\begin{corollary}
  The covariance of $F_n$ and $U_n$ is
  \begin{equation*}
	\cov_n(F_n,U_n)=\E_n(F_nU_n)-\E_nF_n\E_nU_n=-(\Hh_n^{(2)}-1).
  \end{equation*}
\end{corollary}

\begin{proof}
  We have $\cov(J_k,I_k)=-\E J_k\E I_k=-k^{-2}$
  for $k\ge2$, and $\cov(J_1,I_1)=0$, and thus by \refT{T-FU}
  $\cov(F_n,U_n)=\sum_{k=2}^nk^{-2}$.
\end{proof}

It follows further easily from the central limit theorem for
vector-valued random variables, e.g., with Lyapunov's condition, that
the normal limits in Theorems \ref{thm:U} and
\ref{thm:fn} hold jointly, with the joint limit being a pair of two
independent standard normal variables. We omit 
the details.

\section {The number of rows} \label{SR}
In this section we consider  the number
of rows $R_n$  in a random permutation tableau of length $n$. 
Since the rows of a permutation tableau correspond to  \Sx{}
steps in the process of its construction we can write 
$$R_n=\sum_{k=1}^nI_{M_k=\Sx}.$$
\begin{lemma}\label{lem:Smove} For $1\le k\le n$ and any complex
  numbers $z$,  $w$ we have 
 $$\E_k\bigpar{z^{R_k}w^{U_k}}=\frac{w(z-1)}k
\E_{k-1}\bigpar{z^{R_{k-1}}w^{U_{k-1}}}+
\frac{w}k
\E_{k-1}\bigpar{z^{R_{k-1}}(w+1)^{U_{k-1}}}.
$$
\end{lemma}
\begin{proof} By conditioning on $\cF_{k-1}$ we get
\begin{equation}\label{cond4r}
\E_k(z^{R_k}w^{U_k})=
\E_k\left(z^{R_{k-1}}\E_k(z^{I_{M_k=\Sx}}w^{U_k}\mid \cF_{k-1})\right).
\end{equation}
Note that
 $M_k=\Sx$ if and only if $U_k=1+U_{k-1}$ and,  conditionally on $\cF_{k-1}$,
this happens with probability $1/2^{U_{k-1}}$. 
Hence, using also \refL{iter-F} with $z=1$ (or the
 argument in \eqref{gun}),
\begin{align*}
\E_k(z^{I_{M_k=\Sx}}w^{U_k}\mid \cF_{k-1})
&=
\E_k\bigpar{(z^{I_{M_k=\Sx}}-1)w^{U_k}\mid \cF_{k-1}}
+ \E_k(w^{U_k}\mid \cF_{k-1})
\\&=
\E_k\bigpar{(z-1){I_{M_k=\Sx}}w^{U_k}\mid \cF_{k-1}}+
w\left(\frac{w+1}2\right)^{U_{k-1}}
\\&=
(z-1)\frac{w^{1+U_{k-1}}}{2^{U_{k-1}}}+
w\left(\frac{w+1}2\right)^{U_{k-1}}
\\&=
w(z-1)\left(\frac{w}2\right)^{U_{k-1}}+
w\left(\frac{w+1}2\right)^{U_{k-1}}.
\end{align*}
Putting this into (\ref{cond4r}) and applying (\ref{relation}) proves the 
lemma.
\end{proof}

We can now compute the \pgf{} of $R_n$. 
\begin{theorem}\label{T:R}
We have
$$
g_{R_n}(z)
=\frac1{n!}\sum_{r=1}^{n}\euler{n}{r-1} \, z^{r},
$$
where
$\euler nk$ 
are the Eulerian numbers counting the number of permutations
$\sigma=(\sigma_1,\dots,\sigma_n)$ of $[n]$ with $k$ descents (recall
that a descent is a pair $(\sigma_j,\sigma_{j+1})$ such that
$\sigma_j>\sigma_{j+1})$.  
Thus $R_n$ has the same distribution as $1$ plus the
number of descents in a random permutation of $[n]$.
\end{theorem}
\begin{proof}  Consider $\E_nz^{R_n}$.
Applying Lemma~\ref{lem:Smove} $k$ times 
and collecting the terms 
involving the expectation of the same expression we see that it is of
the form 
 \begin{equation}\label{kstep:R}
\E_nz^{R_n}=\frac{1}{n\dotsm(n-k+1)}\sum_{m=0}^kc_{k,m}(z)
\E_{n-k}\bigpar{z^{R_{n-k}}(m+1)^{U_{n-k}}},
\end{equation}
(for $0\le k\le n$)
with certain coefficient functions $c_{k,m}=c_{k,m}(z)$.
Apply Lemma~\ref{lem:Smove}  again
to each of the expectations in the sum  
to get
\begin{align*}
\E_{n-k}\bigpar{z^{R_{n-k}}(m+1)^{U_{n-k}}}
&=
\frac{(m+1)(z-1)}{n-k}\E_{n-k-1}\bigpar{z^{R_{n-k-1}}(m+1)^{U_{n-k-1}}}
 \\&\qquad+\frac{m+1}{n-k}\E_{n-k-1}\bigpar{z^{R_{n-k-1}}(m+2)^{U_{n-k-1}}}.
\end{align*}
Putting that back in, we see that $(n\dotsm(n-k))g_{R_n}(z)$ is equal to
 \begin{align*}
&\sum_{m=0}^kc_{k,m}(m+1)(z-1)\E_{n-k-1}\bigpar{z^{R_{n-k-1}}
(m+1)^{U_{n-k-1}}}
\\&\quad\quad+
\sum_{m=0}^kc_{k,m}(m+1)\E_{n-k-1}\bigpar{z^{R_{n-k-1}}
(m+2)^{U_{n-k-1}}}.
\end{align*}
By rearranging the terms this is
\begin{align*}
&c_{k,0}(z-1)\E_{n-k-1}\bigpar{z^{R_{n-k-1}}}
\\&\qquad\qquad+
\sum_{m=1}^kc_{k,m}(m+1)(z-1)\E_{n-k-1}\bigpar{z^{R_{n-k-1}}
(m+1)^{U_{n-k-1}}}
\\&\quad\quad+
\sum_{m=0}^{k-1}c_{k,m}(m+1)\E_{n-k-1}\bigpar{z^{R_{n-k-1}}
(m+2)^{U_{n-k-1}}}
\\&\qquad\qquad+c_{k,k}(k+1)
\E_{n-k-1}\bigpar{z^{R_{n-k-1}}
(k+2)^{U_{n-k-1}}}\\&\quad
=
c_{k,0}(z-1)\E_{n-k-1}z^{R_{n-k-1}}
\\&\qquad+
\sum_{m=1}^k\left\{c_{k,m}(m+1)(z-1)+mc_{k,m-1}\right\}
\E_{n-k-1}\bigpar{z^{R_{n-k-1}}
(m+1)^{U_{n-k-1}}}
\\&\qquad\qquad+c_{k,k}(k+1)
\E_{n-k-1}\bigpar{z^{R_{n-k-1}}
(k+2)^{U_{n-k-1}}}
\\&\quad
=
\sum_{m=0}^{k+1}c_{k+1,m}\E_{n-k-1}\bigpar{z^{R_{n-k-1}}
(m+1)^{U_{n-k-1}}},
\end{align*}
where the coefficients $c_{k,m}$ satisfy the following recurrence:
$c_{k,m}=0$ unless $0\le m\le k$, $c_{0,0}=1$,  and for $0\le m\le k$, 
$$c_{k,m}=mc_{k-1,m-1}+(m+1)(z-1)c_{k-1,m}.$$
It follows by induction that $c_{k,m}(z)$ is a constant
times $(z-1)^{k-m}$, so 
$$c_{k,m}=a_{k,m}(z-1)^{k-m},$$ 
where $a_{0,0}=1$ and
\begin{equation}\label{a-recur}a_{k,m}=
  \begin{cases}
ma_{k-1,m-1}+(m+1)a_{k-1,m},& \mbox{if\ } 0\le m\le k,\ k\ge1;\\ 
0,&\mbox{otherwise}.	
  \end{cases}
\end{equation}

If we now let   $k=n-1$ in (\ref{kstep:R}) and use
$\E_1\bigpar{z^{R_1}(m+1)^{U_1}}=z(m+1)$ then we get 
\begin{equation}\label{eq:R-a}
g_{R_n}(z)=\frac{z}{n!}\sum_{m=0}^{n-1}(m+1)c_{n-1,m}
=\frac{z}{n!}\sum_{m=0}^{n-1}(m+1)a_{n-1,m}(z-1)^{n-1-m}.
\end{equation}

The recurrence \eqref{a-recur} is solved by
$$a_{k,m}=m!\,\stir{k+1}{m+1},$$
where $\stir{k}m$ is the number of partitions of the $k$ element set
into $m$ non-empty subsets (recall the basic recurrence
$\stir{k+1}{m+1}=(m+1)\stir{k}{m+1}+\stir km$). 
Hence, (\ref{eq:R-a}) becomes 
\begin{align*}g_{R_n}(z)&=\frac{z}{n!}\sum_{m=0}^{n-1}(m+1)!
\,\stir n {m+1}\,(z-1)^{n-(m+1)}\\&= 
\frac{z}{n!}\sum_{m=0}^{n}m!\,\stir nm\,(z-1)^{n-m},
\end{align*}
where in the last step we shifted the index by one and used
$\stir n0=0$ for $n\ge1$. This completes the proof
since for a complex $z$ (see \cite[Section~5.1.3]{knuth}) 
\begin{equation*}
\sum_{r=1}^n\euler n{r-1}\,\frac{z^r}{n!}
=\frac{z}{n!}\sum_{m=0}^{n}m!\,\stir n m(z-1)^{n-m}.
\qedhere
\end{equation*}
\end{proof}

Once the coefficients have been identified as Eulerian numbers we can
use their known properties  
(see e.g.\ \cite[Chapter~10, pp.~150--154]{db}) to obtain
\begin{corollary}\label{cor:anofrn} 
The number of rows $R_n$ in a random permutation tableau of length $n$
satisfies 
\begin{equation}\label{eq:anofrn} 
\frac{R_n-(n+1)/2}{\sqrt{(n+1)/12}}\stackrel d\longrightarrow N(0,1).
\end{equation}
\end{corollary}

In fact,  a local limit theorem  holds as well (see e.g.\ 
\cite{ckss} or \cite{Esseen}).

Moreover, for the number of columns in a random
permutation tableau, by definition, we have $C_n=n-R_n$. Further, let
$D_n$ denote the number of descents in a random permutation;
thus $R_n\eqd D_n+1$ by \refT{T:R}. Hence $C_n=n-R_n\eqd
n-1-D_n=A_n$, the number of ascents in a random permutation. By
symmetry, $A_n\eqd D_n\eqd R_n-1$, and thus we obtain the following
symmetry property:

\begin{corollary}
  The number $C_n$ of columns in a random permutation tableau
  satisfies $C_n\eqd R_n-1$.

In particular,
\begin{equation*}
\frac{C_n-(n-1)/2}{\sqrt{(n+1)/12}}\stackrel d\longrightarrow N(0,1).
\end{equation*}
\end{corollary}

The fact that $R_n-1\eqd C_n\eqd D_n\eqd A_n$ follows
also by the bijections described in \cite{CN}.

 \begin{remark} 
The coefficients in (\ref{eq:R-a}) can be more explicitly written as 
  \begin{equation}\label{coefanr}
 a_{n-1,n-r-1}=(n-r-1)!
 \sum_{1\le j_1<\dots< j_r\le n-1}j_1(j_2-1)\dotsm (j_r-(r-1)).
 \end{equation}
  This can be seen by   putting them in a Pascal type triangle
 $$
 \begin{array}{ccccc}
 & &a_{0,0} &&\\
 &&\swarrow\hspace{.3cm}\searrow&&\\
 &&&&\\
 &a_{1,0}&&a_{1,1}&\\
 &&&&\\
 &\swarrow\hspace{.4cm}\searrow&&\swarrow\hspace{.3cm}\searrow&\\
  &\dots&\dots&\dots&\\
  &\dots&\dots&\dots&\\
 \swarrow &\searrow\hspace{.2cm}\swarrow&\dots&\searrow\hspace{.2cm}\swarrow&\searrow
 \\
 a_{n-1,0} &a_{n-1,1}&\dots&a_{n-1,n-2}&a_{n-1,n-1}
 \end{array}
 $$
 and observing that by (\ref{a-recur}), a move down (either SW or SE)
 from a coefficient $a_{k,m}$ has  weight  $m+1$. The  value of a
 given coefficient at the bottom is obtained by summing over all
 possible paths leading to it from the root $a_{0,0}$ the products of
 weights  corresponding to the moves along the path.   
 Any path contributing to $a_{n-1,n-1-r}$, $0\le r\le n-1$ has exactly
 $r$ SW moves and if they are from  levels $j_1,\dots,j_r$
(with level 0 at the top and $n-1$ at the bottom),
 then the path  has  weight 
 $$(n-1-r)!(j_1+1)j_2(j_3-1)\dotsm(j_r-(r-2)),$$ and (\ref{coefanr}) follows by renumbering the terms. 
 
 Expression (\ref{coefanr}) gives a direct way to get the moments;
 from (\ref{eq:R-a})  the probability generating function of $R_n-1$
 is 
 $$g_{R_n-1}(z)=
\frac1{n!}\sum_{r=0}^{n-1}(n-r)a_{n-1,n-r-1}(z-1)^r,$$
 and thus the 
$r$th factorial moment of $R_n-1$ is 
\begin{align*}&
\E_n\bigpar{(R_n-1)(R_n-2)\dotsm(R_n-r)}
=\frac{d^r}{dz^r}g_{R_n-1}(z)\Big|_{z=1}=\frac{r!(n-r)}{n!}a_{n-1,n-r-1}
\\&\quad=\frac{r!\,(n-r)!}{n!}
 \sum_{1\le j_1<\dots< j_r\le n-1}j_1(j_2-1)\dotsm (j_r-(r-1)).
\end{align*}
In particular, 
$$\E_nR_n=1+\E_n(R_n-1)=1+\frac1n\sum_{j=1}^{n-1}j=1+\frac{n-1}2=\frac{n+1}2,$$
and
\begin{align*}\var(R_n)&=\var(R_n-1)
=\E_n((R_n-1)(R_n-2))+\E_n(R_n-1)-(\E_n(R_n-1))^2\\&=
\frac2{(n-1)n}\sum_{1\le j_1<j_2\le n-1}j_1(j_2-1)+\frac{n-1}2-\frac{(n-1)^2}4\\&=
\frac{(n-2)(3n-5)}{12}+\frac{n-1}2-\frac{(n-1)^2}4=
\frac{n+1}{12},\end{align*}
which explains the normalization in (\ref{eq:anofrn}). 
\end{remark}

\section {Probability generating function of the  number of superfluous 1s}
Let $S_n$ be the number of superfluous 1s in a random tableau of
length $n$.  In \cite{ch} calculations  based on (\ref{relation}) and
(\ref{conddist}) were used to show that  
\begin{equation}\E_nS_n=\frac{(n-1)(n-2)}{12}.\label{expof1s}\end{equation}
In this section we use the same approach to derive the expression for
the \pgf{} 
of $S_n$, see
Proposition~\ref{prop:chfofsn}  or equation (\ref{eq:chfofsn})
below. Although the form of the \pgf{} looks rather
unwieldy,    it can be used  
to compute  (in practice low order) moments of $S_n$ in a relatively
straightforward, albeit a bit tedious way. We will illustrate it by
deriving the exact expression for the variance of  $S_n$ (see
Proposition~\ref{prop:varsn}) below. 
However, there does not seem to be an easy way to derive the central
limit theorem for $S_n$ directly from the \pgf. For
this reason, in the next section we will rely on a different approach
to establish the CLT for $S_n$.

We write in this section
$$S_n=\sum_{k=1}^nV_k,$$ where $V_k$ is the increase in the number of
superfluous 1s when the length of a tableau is increased from $k-1$
to $k$. As in the preceding section let $G_k$ be the row number of the
topmost 1 (provided that the $k$th move is \W{} and counting
only unrestricted rows) and recall that it is set to be $U_{k-1}+1$ if
the $k$th move is \Sx. Then by the same argument as in
Lemma~\ref{lem:G} 
$$\P_k(G_k=j\mid \cF_{k-1})=\frac1{2^j},\quad j=1,\dots,U_{k-1},$$ 
and $$\P_k(G_k=U_{k-1}+1\mid \cF_{k-1})=\frac1{2^{U_{k-1}}}.$$
Moreover, the joint conditional distribution of $(U_k,V_k)$ given
$\cF_{k-1}$ is 
\begin{equation}
\mathcal L((U_k,V_k)\mid \cF_{k-1})
=(G_k+\Bin(U_{k-1}-G_k),\Bin(U_{k-1}-G_k)),\label{jointuv}
\end{equation}
in the sense that if we further condition on $G_k=m$, then the
distribution equals that of $(m+X,X)$ with
$X\sim\Bin(U\kk-m)$ (thus the two occurences of
$\Bin(U\kk-G_k)$ in \eqref{jointuv} signify the
same random variable); further we interpret $\Bin(-1)$ as 0.

\begin{lemma} \label{L5a}
We have, for all complex $z$ and $w$, and
  $k\ge1$:
   \begin{multline*}
\E_k(z^{V_k}w^{U_k}\mid \cF_{k-1})
\\
=w\left\{
\frac1{1+w(z-1)}\left(\frac{zw+1}2\right)^{U_{k-1}}
+\left(1-\frac1{1+w(z-1)}\right)\left(\frac w2\right)^{U_{k-1}}\right\}.  
   \end{multline*}
and
\begin{multline*}
\E_k(z^{S_k}w^{U_k})
=\frac wk\biggl\{
\frac1{1+w(z-1)}\E\kk\left(z^{S\kk}(zw+1)^{U_{k-1}}\right)
\\
+\Bigpar{1-\frac1{1+w(z-1)}}\E\kk\left(z^{S\kk}w^{U_{k-1}}\right)
\biggr\}.   
\end{multline*}
\end{lemma}
\begin{proof} Using (\ref{jointuv}) we have
\begin{align*}
&\E_k(z^{V_k}w^{U_k}\mid \cF_{k-1})=\E_k( z^{\Bin(U_{k-1}-G_k)}w^{G_k+\Bin(U_{k-1}-G_k)}\mid \cF_{k-1})\\
&\qquad=\sum_{m=1}^{U_{k-1}}
\E_k\bigpar{(zw)^{\Bin(U_{k-1}-m)}w^{m}I(G_k=m)\mid
  \cF_{k-1}}\\&\qquad\qquad\qquad+w^{1+U_{k-1}}\P_k(G_k= U_{k-1}+1\mid
\cF_{k-1})\\ 
&\qquad=\sum_{m=1}^{U_{k-1}}\left(\frac w2\right)^m\E_k((zw)^{\Bin(U_{k-1}-m)}\mid \cF_{k-1})+w\left(\frac w2\right)^{U_{k-1}}
 \\
 &\qquad=\sum_{m=1}^{U_{k-1}}\left(\frac w2\right)^m\left(\frac{zw+1}2\right)^{U_{k-1}-m}+w\left(\frac w2\right)^{U_{k-1}}\\
 &\qquad=\sum_{\ell=0}^{U_{k-1}-1}\left(\frac{zw+1}2\right)^{\ell}\left(\frac w2\right)^{U_{k-1}-\ell}
 +w\left(\frac w2\right)^{U_{k-1}}\\
&\qquad=\left(\frac w2\right)^{U_{k-1}}\left\{
\frac 1{(1+zw)/w-1}
\left(\left(\frac{zw+1}w\right)^{U_{k-1}}-1\right)+w\right\}\\
&\qquad=w\left\{
\frac1{1+w(z-1)}\left(\frac{zw+1}2\right)^{U_{k-1}}
+\left(1-\frac1{1+w(z-1)}\right)\left(\frac w2\right)^{U_{k-1}}\right\},
\end{align*}
which is the first formula.

To prove the second write
\begin{align*}\E_{k}\bigpar{z^{S_{k}}w^{U_{k}}}&=\E_{k}\E_{k}(
z^{S_{k}}w^{U_{k}}\mid \cF_{k-1})
=\E_{k}\bigpar{z^{S_{k-1}}
\E_{k}(z^{V_{k}}w^{U_{k}}\mid \cF_{k-1})}
\end{align*}
and use the first part and the usual reduction by
\eqref{relation}
from $\P_{k}$ to $\P_{k-1}$.
\end{proof}

For $\ell\ge0$ set
$b_\ell=b_\ell(z)=\sum_{j=0}^\ell z^{j}$ so that we have
 \begin{equation}\label{b}zb_\ell+1=b_{\ell+1}\quad\mbox{and}\quad
   1+(z-1)b_\ell=z^{\ell+1}. 
\end{equation}
If we substitute $b_\ell$ for $w$ in \refL{L5a} and use
\eqref{b} we obtain the following.
\begin{lemma} 
For $0\le k\le n-1$ and $\ell\ge0$,
\begin{multline*}
\E_{n-k}\bigpar{z^{S_{n-k}}b_\ell^{U_{n-k}}}
=\frac{b_\ell}{n-k}\left(z^{-\ell-1}
\E_{n-k-1}\bigpar{z^{S_{n-k-1}}b_{\ell+1}^{U_{n-k-1}}}\right.
\\
\left. {}+(1-z^{-\ell-1})\E_{n-k-1}\bigpar{z^{S_{n-k-1}}b_{\ell}^{U_{n-k-1}}}\right).
\end{multline*}
\end{lemma}

\begin{proposition}\label{prop:chfofsn}
We have
$$g_{S_n}(z)=\frac1{n!}\sum_{\ell=0}^{n-1}c_{n-1,\ell}(z)b_{\ell}(z),$$
where the coefficients $c_{m,\ell}=c_{m,\ell}(z)$ satisfy 
$c_{0,0}=1$
and for $m\ge0$ the recurrence 
$$c_{m+1,\ell}
=c_{m,\ell}(1-z^{-\ell-1})b_\ell(z)+c_{m,\ell-1}z^{-\ell}b_{\ell-1}(z).$$
In particular, $c_{m,\ell}=0$ unless
$0\le\ell\le m$.
\end{proposition}

\begin{proof}  Skipping the dependence of $b_\ell$'s on $z$ for
  notational convenience, we show first that for $0\le k\le n-1$ we have
\begin{equation}
\label{kstep-S}\E_nz^{S_n}=\E_nz^{S_n}b_0^{U_n}=\frac1{n\dotsm(n-k+1)}\sum_{m=0}^kc_{k,m}\E_{n-k}z^{S_{n-k}}b_m^{U_{n-k}}.
\end{equation}
Indeed, assume inductively \eqref{kstep-S} for some
$k\ge0$ ($k=0$ being trivial).
Apply the previous lemma to each of the terms $\E_{n-k}z^{S_{n-k}}b_m^{U_{n-k}}$ to get
\begin{align*}\E_nz^{S_n}&=\frac1{n\dotsm(n-k)}\sum_{m=0}^kc_{k,m}b_m\left\{z^{-m-1}\E_{n-k-1}z^{S_{n-k-1}}b_{m+1}^{U_{n-k-1}}\right.\\ &\qquad \qquad \qquad \left. +
(1-z^{-m-1})\E_{n-k-1}z^{S_{n-k-1}}b_{m}^{U_{n-k-1}}\right\}
.\end{align*}
Separating the sums, rearranging the terms, and collecting the coefficients in front of $\E_{n-k-1}z^{S_{n-k-1}}b_{m}^{U_{n-k-1}}$ we obtain that $\E_nz^{S_n}$ is equal to 
\begin{align*}
& \frac1{n\dotsm(n-k)}\Big\{
\sum_{m=0}^{k+1}\left(c_{k,m-1}b_{m-1}z^{-m}+c_{k,m}b_m(1-z^{-m-1})\right)\E_{n-k-1}z^{S_{n-k-1}}b_{m}^{U_{n-k-1}}\Big\}
.\end{align*}
This completes the inductive proof of \eqref{kstep-S} with
the coefficients $c_{k,m}$ given by the specified recurrence. 
Choosing $k=n-1$ in (\ref{kstep-S}) and using the observation that
$$\E_1z^{S_1}b_m^{U_1}=b_m,$$
(because $S_1=0$ and $U_1=1$) completes the proof.
\end{proof}

 \begin{remark} Just as in the proof of  \refT{T:R}, 
the coefficients  $c_{k,m}$ can be 
 put in a Pascal-type triangle
 $$
 \begin{array}{ccccc}
 & &c_{0,0} &&\\
 &&\swarrow\hspace{.3cm}\searrow&&\\
 &&&&\\
 &c_{1,0}&&c_{1,1}&\\
 &&&&\\
 &\swarrow\hspace{.4cm}\searrow&&\swarrow\hspace{.3cm}\searrow&\\
  &\dots&\dots&\dots&\\
  &\dots&\dots&\dots&\\
 \swarrow &\searrow\hspace{.2cm}\swarrow&\dots&\searrow\hspace{.2cm}\swarrow&\searrow
 \\
 c_{n-1,0} &c_{n-1,1}&\dots&c_{n-1,n-2}&c_{n-1,n-1}
 \end{array}
 $$
 This time
a SE move from a coefficient $c_{m,\ell}$ has  weight 
$z^{-\ell-1}b_\ell(z)$ and a SW move  has  weight $(1-z^{-\ell-1})b_\ell(z)$. 
As was in the case of the number of rows, the  value of a given
 coefficient at the bottom is obtained by summing over all possible
 paths leading to it from the root $c_{0,0}$ the products of weights
 corresponding to the moves along the path.  This can be used to
 obtain an explicit expressions for the coefficients $(c_{n-1,m})$.
 Any path from $c_{0,0}$ to $c_{n-1,n-1-r}$, $0\le r\le n-1$, has
 exactly $r$ SW moves and $n-1-r$  
 SE moves. The total weight of SE moves (they are from $c_{\cdot,0},c_{\cdot,1},\dots,c_{\cdot,n-2-r}$) is
 \[\prod_{j=0}^{n-2-r}z^{-j-1}b_j(z)
=z^{-\binom{n-r}2}\prod_{j=0}^{n-r-2}b_j(z).\]
 The SW moves may be from  $c_{\cdot,\ell_1},\dots,c_{\cdot,\ell_r}$, for some $0\le\ell_1\le\dots\le\ell_r\le n-1-r$. A SW move from $c_{\cdot,\ell_j}$ has weight $(1-z^{-\ell_j-1})b_{\ell_j}(z)$. Therefore,
 \[c_{n-1,n-1-r}(z)=z^{-\binom{n-r}2}\prod_{k=0}^{n-r-2}b_k(z)\sum_{0\le\ell_1\le\dots\le\ell_r\le n-1-r}\prod_{j=1}^r(1-z^{-\ell_j-1})b_{\ell_j}(z).
 \]
This gives an explicit expression for the \pgf{} of $S_n$, namely
\begin{equation}\label{eq:chfofsn}
g_{S_n}(z)=
\frac1{n!}\sum_{r=0}^{n-1}z^{-\binom{n-r}2}
\left(\prod_{k=0}^{n-r-1}b_k(z)\right)A_{n,r}(z), 
\end{equation}
where 
\begin{equation}\label{eq:anr}A_{n,r}(z)=\sum_{0\le\ell_1\le\dots\le\ell_r\le n-1-r}\prod_{j=1}^r(1-z^{-\ell_j-1})b_{\ell_j}(z).
 \end{equation}
 \end{remark}
(Note that all negative powers cancel in \eqref{eq:chfofsn}
 since $g_{S_n}(z)$ is a polynomial in $z$.)
 Although this expression looks quite complicated and it is not clear
 to us at the moment how to deduce the asymptotic normality of $S_n$
 from it, some information can be extracted from it.  We will
 illustrate it by deriving an exact expression for the variance of
 $S_n$.  
 \begin{proposition}\label{prop:varsn} For $n\ge2$, the
 variance of $S_n$ satisfies 
 \begin{equation}\label{eq:varsn}\var(S_n)=\frac{(n-2)(2n^2+11n-1)}{360}.\end{equation}
 \end{proposition}
 \begin{proof}
We  compute the second factorial moment of $S_n$:
\begin{equation}\label{eq:2ndfac}\E_nS_n(S_n-1)=\frac1{n!}\frac
  {d^2}{dz^2}\Big(\sum_{k=0}^{n-1}c_{n-1,k}b_{k}\Big)\Big|_{z=1}.
\end{equation} 
Notice that every path  contributing to $c_{n-1,r}$, for $0\le r<n-3$ has at least three SW moves and so its weight will have at least three factors of the form $1-z^{-j}$. Hence if it is differentiated twice and evaluated at $z=1$ it will vanish. It follows that the sum on the right hand side of (\ref{eq:2ndfac}) reduces to the last three terms.
For $k=n-3$ we have
$$\frac1{n!}c_{n-1,n-3}b_{n-3}=\frac{z^{-\binom{n-2}2}}{n!}\left(\prod_{k=0}^{n-4}b_k\right)b_{n-3}A_{n,2}(z)=G(z)A_{n,2}(z),$$
where we have set
$$G(z):=\frac1{n!}z^{-\binom{n-2}2}\prod_{k=0}^{n-3}b_{k}(z).$$ 
So, the second derivative of $G(z)A_{n,2}(z)$ is
$$G''(z)A_{n,2}(z)+2G'(z)A'_{n,2}(z)+G(z)A''_{n,2}(z).$$
Since
$$A_{n,2}(z)=\sum_{1\le\ell\le m\le n-2}(1-z^{-\ell})(1-z^{-m})b_{\ell-1}(z)b_{m-1}(z),$$
$A_{n,2}(1)=A'_{n,2}(1)=0$ so we only need $G(1)A''_{n,2}(1)$. 
Now, 
\begin{equation}\label{g(1)}G(1)=\frac1{n!}\prod_{k=0}^{n-3}(k+1)
=\frac1{n(n-1)}. 
\end{equation}
To compute $A''_{n,2}(1)$, writing 
$$(1-z^{-\ell})(1-z^{-m})=h_1(z),\quad b_{\ell-1}(z)b_{m-1}(z)=h_2(z),$$
we see that
$$h''_1(1)h_2(1)+2h'_1(1)h'_2(1)+h_1(1)h''_2(1)=h''_1(1)h_2(1)=2\ell^2m^2.$$
Therefore,
\begin{align*}A''_{n,2}(1)&=2\sum_{m=1}^{n-2}m^2\sum_{\ell=1}^{m}\ell^2
=\frac13\sum_{m=1}^{n-2}m^3(m+1)(2m+1)\\&
=\frac1{180}n(n-1)(n-2)(5n-11)(2n-1)(2n-3),
\end{align*}
and combining this with (\ref{g(1)}) we obtain
$$\frac1{n!}\frac{d^2}{dz^2}\left(c_{n-1,n-3}b_{n-3}\right)\Big|_{z=1}
=\frac{(n-2)(5n-11)(2n-1)(2n-3)}{180}
.$$
We next handle 
\begin{align*}\frac1{n!}c_{n-1,n-1}b_{n-1}
&=\frac{z^{-\binom{n}2}}{n!}\left(\prod_{k=0}^{n-2}b_k\right)b_{n-1}
\\&
=z^{-\binom{n}2}\prod_{j=1}^{n-1}\frac{b_{j-1}}j.
\end{align*}
The last expression is the probability generating function of a sum of
 independent random variables $W_0,\dots,W_{n}$, where 
$W_{0}\equiv-n(n-1)/2$ and
for $1\le j\le n$, $W_j$  
 is a discrete uniform random variable on $\{0,\dots,j-1\}$.
So, with $Y_n=\sum_{j=0}^{n}W_j$ we have
$$\frac1{n!}\frac{d^2}{dz^2}c_{n-1,n-1}b_{n-1}\Big|_{z=1}=\E Y_n(Y_n-1)=
\sum_{j=1}^n\var(W_j)+\left(\E Y_n\right)^2-\E Y_n.$$
Now,
$$\E W_j=\frac{j-1}2,\qquad1\le j\le n,
$$
so that
$$\E Y_n=\frac12\sum_{j=1}^{n}(j-1)-\frac{n(n-1)}2=
-\frac{n(n-1)}{4}.$$
Furthermore,
 for $1\le j\le n$
$$\var(W_j)=\frac1j\sum_{k=0}^{j-1}k^2-\left(\frac{j-1}2\right)^2=
\frac{(j-1)(j+1)}{12}.$$
 Therefore,
 \begin{align*}\sum_{j=1}^{n}\var(W_j)&=\frac{n(n-1)(2n+5)}{72}.
 \end{align*}
Finally, putting the above together we get
\begin{align*}&\frac1{n!}\frac{d^2}{dz^2}c_{n-1,n-1}b_{n-1}\Big|_{z=1}
=
\frac{n(n-1)(2n+5)}{72}
+\frac{n^2(n-1)^2}{16}
+\frac{n(n-1)}{4}
\\&\quad=
\frac{n(n-1)(9n^2-5n+46)}{144}.
\end{align*}
It remains to handle $\frac1{n!}c_{n-1,n-2}b_{n-2}$. 
 We write
$$c_{n-1,n-2}b_{n-2}
=
f(z)A_{n,1}(z),
$$
where
\[f(z)=z^{-\binom{n-1}2}\prod_{k=0}^{n-2}b_k(z)\]
and, according to (\ref{eq:anr}),
\[A_{n,1}(z) =\sum_{\ell=0}^{n-2}(1-z^{-\ell-1})b_{\ell}(z).\]
Then the second derivative of $c_{n-1,n-2}b_{n-2}$ is
$$f''(z)A_{n,1}(z)+2f'(z)A'_{n,1}(z)+f(z)A''_{n,1}(z).$$
At $z=1$ the first product vanishes. For the remaining two first notice that  
$$\frac1{(n-1)!}f(z)
=z^{-\binom{n-1}2}\prod_{j=1}^{n-1}\frac{b_{j-1}}j,
$$
is a generating function of a legitimate distribution function. Therefore,
$f(1)=(n-1)!$ and 
 $f'(1)/(n-1)!$ is the expected value of a random variable represented by $f(z)/(n-1)!$. Since this expected value  is
$$-\binom{n-1}2+\sum_{j=1}^{n-1}\frac1j\sum_{k=0}^{j-1}k=
-\frac{(n-1)(n-2)}4,$$
we obtain
$$f'(1)=-\frac{(n-1)(n-2)(n-1)!}4.$$
It remains to compute  the first two derivatives of $A_{n,1}(z)$ at $z=1$.
$$A'_{n,1}(z)=\sum_{\ell=1}^{n-1}\left(b'_{\ell-1}(1-z^{-\ell})+b_{\ell-1}\ell z^{-\ell-1}\right),$$
so that
\begin{equation}\label{fprime2}
A'_{n,1}(1)=\sum_{\ell=1}^{n-1}\ell^2=\frac{n(n-1)(2n-1)}6.
\end{equation}
Also
$$A''_{n,1}(z)=\sum_{\ell=1}^{n-1}\left(b''_{\ell-1}(1-z^{-\ell})+2b'_{\ell-1}\ell z^{-\ell-1}-b_{\ell-1}\ell(\ell+1)z^{-\ell-2}\right).$$
Since $b'_{\ell-1}(1)=(\ell-1)\ell/2$  we get   
$$A''_{n,1}(1)=\sum_{\ell=1}^{n-1}\left((\ell-1)\ell^2-\ell^2(\ell+1)\right)
=-\frac{n(n-1)(2n-1)}3.$$ 
Hence,
\begin{align*}
\frac1{n!}\frac{d^2}{dz^2}c_{n-1,n-2}b_{n-2}\Big|_{z=1}
&=\frac1{n!}\bigpar{2f'(1)A'_{n,1}(1)+f(1)A''_{n,1}(1)}
\\&
=-\frac{(n-1)^2(n-2)(2n-1)}{12}
-\frac{(n-1)(2n-1)}{3}
\\&=
-\frac{(n-1)(2n-1)(n^2-3n+6)}{12}.
\end{align*}
Combining all of these calculations  gives 
\begin{align*}
\E_nS_n(S_n-1)
&=
\frac{(n-2)(5n-11)(2n-1)(2n-3)}{180}
\\&\qquad
-\frac{(n-1)(2n-1)(n^2-3n+6)}{12}
\\&\qquad
+\frac{n(n-1)(9n^2-5n+46)}{144}
\\&=
\frac{(n-2)(n-3)(5n^2-n-16)}{720}.
\end{align*}
By the same argument
\begin{align*}\E_nS_n
&=\frac1{n!}\frac d{dz}
\Big(c_{n-1,n-1}b_{n-1}+c_{n-1,n-2}b_{n-2}\Big)\Big|_{z=1}
\\&=\E Y_n+
\frac1{n!}\left(f'(1)A_{n,1}(1)+f(1)A'_{n,1}(1)\right)
\end{align*}
Using $A_{n,1}(1)=0$, $f(1)=(n-1)!$, (\ref{fprime2}), and the value of
$\E Y_n$ we get  
$$
\E_nS_n=-\frac{n(n-1)}{4}+\frac{(n-1)(2n-1)}{6}=\frac{(n-1)(n-2)}{12},$$
which conforms to (\ref{expof1s}).
Finally,
\begin{align*}\var(S_n)&=\E S_n(S_n-1)-(\E S_n)^2+\E S_n\\&=
\frac{(n-2)(n-3)(5n^2-n-16)}{720}-\left(\frac{(n-1)(n-2)}{12}\right)^2+
\frac{(n-1)(n-2)}{12}\\&=
\frac{(n-2)(2n^2+11n-1)}{360},
\end{align*}
which proves Proposition~\ref{prop:varsn}.
\end{proof}

\
\section{Asymptotic normality of $S_n$} 

In this section we  provide a self-contained proof of  the following
\begin{theorem} \label{TSn}
As $n\to\infty$ we have
\begin{equation}\label{clt41s}
\frac{S_n-n^2/12}{\sqrt{n^3/180}}\stackrel d\longrightarrow N(0,1).
\end{equation}
\end{theorem}
\begin{proof} As we mentioned, the form of the \pgf{} of
  $S_n$ obtained in the previous section  does not seem to be
  convenient to yield the central limit theorem. For this reason,  we
  will rely on a bijective result of Steingr\'{\i}msson and Williams
  \cite{SW}. According to their result the number of superfluous 1s
  in a permutation tableau of length $n$ is equidistributed with the
  number of occurrences of a generalized patterns $31{-}2$ in a random
  permutation  of $[n]$. (An occurrence of a generalized pattern $31{-}2$
  in a permutation $\sigma$ is a pair $1<i<j$ such that
  $\sigma_{i-1}>\sigma_j>\sigma_i$.) To analyze that quantity, it will
  be convenient to think of a random permutation as generated from a
  sample $X_1,\dots, X_n$  i.i.d.\ random variables  
each being uniform on $[0,1]$ (the permutation is obtained by reading
  off the ranks of $X_1,\dots,X_n$). If for $2\le i<j\le n$ we let
  $I_{i,j}:=I_{X_{i-1}>X_j>X_{i}}$ then  
$$S_n=\sum_{2\le i<j\le n}I_{i,j}.$$ 
Notice that from this representation we immediately recover
(\ref{expof1s}) since 
$$\E S_n=\E\sum_{2\le i<j\le n}I_{i,j}=\binom{n-1}2\P(X_1>X_3>X_2)
=\frac16\binom{n-1}2.$$
Similarly, we can easily obtain the asymptotic value of the variance: we write
$$
\var(S_n)=\sum_{\substack{i_1<j_1\\i_2<j_2}}\cov(I_{i_1,j_1},I_{i_2,j_2}),$$
and note that if $\{i_1-1,i_1,j_1\}\cap\{i_2-1,i_2,j_2\}=\emptyset$ then 
$I_{i_1,j_1}$ and $I_{i_2,j_2}$ are independent and so their
covariance vanishes. In the complementary case, the main contribution
comes from the cases that contribute  $\Theta(n^3)$ terms  to the
sum. We obtain 
\begin{multline*}
\var(S_n)\sim\frac{n^3}3\Big(\cov(I_{2,3},I_{2,4})+\cov(I_{2,5},I_{4,5})
+\cov(I_{2,4},I_{3,5})\\
{}+\cov(I_{2,5},I_{3,4})
+
\cov(I_{2,3},I_{4,5})+\cov(I_{2,4},I_{4,5})\Big),
\end{multline*}
as all other cases contribute $O(n^2)$ terms to the sum. We calculate:
$$\E I_{2,3}\cap I_{2,4}=\P(X_1>X_3>X_2,X_1>X_4>X_2)=2\P(X_1>X_3>X_4>X_2)=
\frac2{4!},$$
so that
$$\cov(I_{2,3},I_{2,4})=\frac1{12}-\left(\frac1{6}\right)^2.$$
Similarly,
\begin{align*}
\E I_{2,5}\cap I_{4,5}&=\P(X_1>X_5>X_2,X_3>X_5>X_4)\\&=4\P(X_1>X_3>X_5>X_2>X_4)=\frac1{30},\\
\E I_{2,4}\cap I_{3,5}&=\P(X_1>X_4>X_2,X_2>X_5>X_3)=\frac1{120},\\
\E I_{2,5}\cap I_{3,4}&=\P(X_1>X_5>X_2,X_2>X_4>X_3)=\frac1{120},\\
\E I_{2,3}\cap I_{4,5}&=\P(X_1>X_3>X_2,X_3>X_5>X_4)\\&=3\P(X_1>X_3>X_2>X_5>X_4)=\frac1{40},\\
\E I_{2,4}\cap I_{4,5}&=\P(X_1>X_4>X_2,X_3>X_5>X_4)=\frac1{40}.
\end{align*}
Hence,
\begin{equation}\label{varof1s-as}
\var(S_n)\sim\frac{n^3}3\cdot\Bigpar{\frac1{12}+\frac1{30}
+\frac1{120}+\frac1{120}+\frac1{40}+\frac1{40}-\frac6{36}}
=\frac{n^3}{180}.\end{equation}

\begin{remark} 
The exact value of the variance could be obtained by computing the
other terms.  
\end{remark}

Finally, to establish (\ref{clt41s}) we will rely on results presented
in \cite{SJ58} and \cite[Section~6.1]{jlr}; see \cite{Bona} for a
closely related theorem proved by the same method. See also \cite{Esseen} 
for a related simple proof of the asymptotic normality of the number
of descents, cf.\ \refC{cor:anofrn} above.
We let
$A=A_n:=\{\a=(i,j):\ 2\le i<j\le 
n\}$. Then $S_n=\sum_{\a\in A}I_\a$. Recall that a dependency graph
$L$ for $\{I_\a:\ \a\in A\}$ is any graph whose vertex set is $A$ and
which has the property that if $V_1$, $V_2$ are  two disjoint subsets
of $A$ such that $L$ has no edges with one endpoint in $V_1$ and the
other in $V_2$ then the families $\{I_\a:\ \a\in V_1\}$ and $\{I_\a:\
\a\in V_2\}$ are mutually independent. For our purposes it is enough
to consider $L$ defined by the following rule: we put an edge between
$\a_1=(i_1,j_1)$ and $\a_2=(i_2,j_2)$ iff
$\{i_1-1,i_1,j_1\}\cap\{i_2-1,i_2,j_2\}\ne\emptyset$. If
$\a_1,\dots,\a_r\in A$ then the closed neighborhood of
$\{\a_1,\dots,\a_r\}$ in $L$ is defined by 
$$\bar{N}_L(\a_1,\dots,\a_r)=\bigcup_{i=1}^r\{\b\in A:\ \b=\a_i\
\mbox{or}\ \a_i\b\in E(L)\},
$$
where $E(L)$ denotes the edge set of $L$. Note that in our case for
every fixed $r\ge2$ $|\bar{N}_L(\a_1,\dots,\a_{r-1})|=O_r(n)$, where
$O_r(\ \cdot\ )$ means that the constant may depend on $r$. Hence,
trivially 
$$\sum_{\a\in
\bar{N}_L(\a_1,\dots,\a_{r-1})}\E\bigpar{I_\a\big|I_{\a_1},\dots,I_{\a_{r-1}}}
=O_r(n).
$$ 
Since, cf.\ (\ref{expof1s}) for the exact value,
$$\sum_{\a\in A}\E I_\a\le|A|=O(n^2),$$
by \cite[Lemma~6.17]{jlr} we conclude that the $r$th cumulant of
 $S_n$ (defined by
$\kappa_r(S_n)=i^{-k}\frac{d^k}{dt^k}\log\phi_{S_n}(0)$, where  $\phi_X(t)=g_X(e^{it})$ is the characteristic function) satisfies 
$$|\kappa_r(S_n)|=O_r(n^2\cdot n^{r-1})=O_r(n^{r+1}).
$$
Hence, for $r\ge3$ 
$$\left|\kappa_r\left(\frac{S_n-n^2/12}{\sqrt{n^3/180}}\right)\right|
=O_r(n^{r+1}n^{-\frac32r})=o(1),$$
as $n\to\infty$. Since $\kappa_1(X)=\E X$ and $\kappa_2(X)=\var(X)$ we have
$$\kappa_1\left(\frac{S_n-n^2/12}{\sqrt{n^3/180}}\right)\longto0,\quad\quad
\kappa_2\left(\frac{S_n-n^2/12}{\sqrt{n^3/180}}\right)\longto1,$$
and the theorem follows by the cumulant convergence theorem (see
e.g.\ \cite[Corollary~6.15]{jlr}). 
\end{proof}

\begin{remark} Our results can be used to draw conclusions about some
other parameters. For example,  
let $Y_n$ be the number of 1s in the random permutation tableaux
$T_n$.
Although we  have not computed an explicit formula for the
distribution of $Y_n$, we can easily obtain its asymptotic
distribution. Namely, 
as $n\to\infty$ we have
\begin{equation*}
\frac{Y_n-n^2/12}{\sqrt{n^3/180}}\stackrel d\longrightarrow N(0,1).
\end{equation*}
This follows immediately  from \refT{TSn} upon noting that   
\[Y_n=S_n+C_n=S_n+O(n).\]

It seems straightforward to prove a similar central limit theorem for 
$Z_n$, the number of 0s in the random permutation tableaux
$T_n$, using further bijective results by \cite{SW}; the main
difference is that we need to consider several generalized patterns
simultaneously and the joint distribution of their numbers of occurrences. 
We leave this to the reader.
\end{remark}
\bibliographystyle{amsalpha}

\end{document}